\documentclass[11pt]{article}
\usepackage{epsf}

\begin{document}


\newtheorem{theorem}{Theorem}[section]
\newtheorem{itlemma}{Lemma}[section]
\newtheorem{itdefinition}{Definition}[section]
\newtheorem{itexample}{Example}
\newtheorem{itclaim}{Claim}[section]
\newtheorem{itproposition}{Proposition}[section]
\newtheorem{itremark}{Remark}[section]
\newtheorem{itcorollary}{Corollary}[section]

\newenvironment{example}{\begin{itexample}\rm}{\end{itexample}}
\newenvironment{definition}{\begin{itdefinition}\rm}{\end{itdefinition}}
\newenvironment{lemma}{\begin{itlemma}\rm}{\end{itlemma}}
\newenvironment{corollary}{\begin{itcorollary}\rm}{\end{itcorollary}}
\newenvironment{claim}{\begin{itclaim}\rm}{\end{itclaim}}
\newenvironment{proposition}{\begin{itproposition}\rm}{\end{itproposition}}
\newenvironment{remark}{\begin{itremark}\rm}{\end{itremark}}

\newcommand{\qed}{\hfill \halmos} 
\newcommand{\mybox}{\hfill $\Box$} 

\newcommand{\comment}[1]{}
\newcommand{\halmos}{\rule{1ex}{1.4ex}}
\newenvironment{proof}{\noindent {\em Proof}.\ }{\hspace*{\fill}$\halmos$ 
\medskip}

\def\vbar{\mathchoice{\vrule height6.3ptdepth-.5ptwidth.8pt\kern-.8pt}
   {\vrule height6.3ptdepth-.5ptwidth.8pt\kern-.8pt}
   {\vrule height4.1ptdepth-.35ptwidth.6pt\kern-.6pt}
   {\vrule height3.1ptdepth-.25ptwidth.5pt\kern-.5pt}}
\def\fudge{\mathchoice{}{}{\mkern.5mu}{\mkern.8mu}}
\def\bbc#1#2{{\rm \mkern#2mu\vbar\mkern-#2mu#1}}
\def\bbb#1{{\rm I\mkern-3.5mu #1}}
\def\bba#1#2{{\rm #1\mkern-#2mu\fudge #1}}
\def\bb#1{{\count4=`#1 \advance\count4by-64 \ifcase\count4\or\bba A{11.5}\or
   \bbb B\or\bbc C{5}\or\bbb D\or\bbb E\or\bbb F \or\bbc G{5}\or\bbb H\or
   \bbb I\or\bbc J{3}\or\bbb K\or\bbb L \or\bbb M\or\bbb N\or\bbc O{5} \or
   \bbb P\or\bbc Q{5}\or\bbb R\or\bbc S{4.2}\or\bba T{10.5}\or\bbc U{5}\or
   \bba V{12}\or\bba W{16.5}\or\bba X{11}\or\bba Y{11.7}\or\bba Z{7.5}\fi}}

\def\Q{{\bb Q}}                         
\def\N{{\bb N}}                         
\def\R{{\bb R}}                         
\def\I{{\bb Z}}                         
\def\B{{\bb B}}                         

\def\rmtr{{\rm tr}}

\def\nasymptotic{{_{\stackrel{\displaystyle\longrightarrow}
{N\rightarrow\infty}}\,\, }} 
\def\masymptotic{{_{\stackrel{\displaystyle\longrightarrow}
{M\rightarrow\infty}}\,\, }} 
\def\wasymptotic{{_{\stackrel{\displaystyle\longrightarrow}
{w\rightarrow\infty}}\,\, }} 

\def\asymptext{\raisebox{.6ex}{${_{\stackrel{\displaystyle\longrightarro
w}{x\rightarrow\pm\infty}}\,\, }$}} 
\def\epsilim{{_{\textstyle{\rm lim}}\atop_{\epsilon\rightarrow 0+}\,\, }} 

\def\beqra{\begin{eqnarray}} \def\eeqra{\end{eqnarray}}
\def\beqast{\begin{eqnarray*}} \def\eeqast{\end{eqnarray*}}
\def\beq{\begin{equation}}      \def\eeq{\end{equation}}
\def\be{\begin{enumerate}}   \def\ee{\end{enumerate}}

\def\bet{\beta}
\def\gam{\gamma}
\def\Gam{\Gamma}
\def\la{\lambda}
\def\eps{\epsilon}
\def\La{\Lambda}
\def\si{\sigma}
\def\Si{\Sigma}
\def\al{\alpha}
\def\Tha{\Theta}
\def\tha{\theta}
\def\vphi{\varphi}
\def\del{\delta}
\def\Del{\Delta}
\def\ab{\alpha\beta}
\def\om{\omega}
\def\Om{\Omega}
\def\mn{\mu\nu}
\def\mun{^{\mu}{}_{\nu}}
\def\kap{\kappa}
\def\rsi{\rho\sigma}
\def\beal{\beta\alpha}

\def\til{\tilde}
\def\rta{\rightarrow}
\def\eqv{\equiv}
\def\nab{\nabla}
\def\pa{\partial}
\def\sit{\tilde\sigma}
\def\ul{\underline}
\def\indt{\parindent2.5em}
\def\nd{\noindent}
\def\var{{1\over 2\si^2}}
\def\ivar{\left(2\pi\si^2\right)}
\def\iivar{\left({1\over 2\pi\si^2}\right)}
\def\caa{{\cal A}}
\def\cb{{\cal B}}
\def\cac{{\cal C}}
\def\cd{{\cal D}}
\def\ce{{\cal E}}
\def\cf{{\cal F}}
\def\cg{{\cal G}}
\def\ch{{\cal H}}
\def\ci{{\cal I}}
\def\cj{{\cal{J}}}
\def\ck{{\cal K}}
\def\cl{{\cal L}}
\def\cm{{\cal M}}
\def\cn{{\cal N}}
\def\cO{{\cal O}}
\def\cp{{\cal P}}
\def\cq{{\cal Q}}
\def\car{{\cal R}}
\def\cs{{\cal S}}
\def\ct{{\cal{T}}}
\def\cu{{\cal{U}}}
\def\cv{{\cal{V}}}
\def\cw{{\cal{W}}}
\def\cx{{\cal{X}}}
\def\cy{{\cal{Y}}}
\def\cz{{\cal{Z}}}

\renewcommand{\thesection}{\arabic{section}}
\renewcommand{\theequation}{\thesection.\arabic{equation}}

\vspace*{.2in}
\begin{center}
\Large{\sc 
On the universality of the probability distribution of 
the product $B^{-1}X$ of random matrices}\\
\normalsize
\vspace{15pt}
\begin{center}
{\bf Joshua Feinberg $^*$}
\end{center}
\vskip 2mm
\begin{center}
{Physics Department$^{\#}$,}\\
{University of Haifa at Oranim, Tivon 36006, Israel}\\
{and}\\
{Physics Department,}\\
{Technion, Israel Institute of Technology, Haifa 32000, Israel}\\
\vskip 2mm
\end{center}
\vspace{.3cm}
\end{center}

\date{16 March 2002}
\begin{minipage}{5.3in}
{\abstract~~~~~
Consider random matrices $A$, of dimension $m\times (m+n)$,  drawn from
an ensemble with probability density $f(\rmtr AA^\dagger)$, with $f(x)$ 
a given appropriate function. Break $A = (B,X)$ into an $m\times m$ block 
$B$ and the complementary $m\times n$ block $X$, and define the random 
matrix $Z=B^{-1}X$.  We calculate the probability density function 
$P(Z)$ of the random matrix $Z$ and find that it is a universal function, 
independent of $f(x)$. The universal probability distribution $P(Z)$ is a 
spherically symmetric matrix-variate $t$-distribution. Universality of $P(Z)$ 
is, essentially, a consequence of rotational invariance of the probability 
ensembles we study. As an application, we study the distribution 
of solutions of systems of linear equations with random coefficients, and 
extend a classic result due to Girko.} 
\end{minipage}

\vfill
\vspace{10pt}
\flushleft *{\it e-mail address: joshua@physics.technion.ac.il}\\
$\#$ {\it permanent address}\\
AMS {\em subject classifications}: 15A52, 60E05, 62H10, 34F05 
\vfill
\pagebreak

\setcounter{page}{1}

\section{Introduction}
In this note we will address the issue of universality of the 
probability density function (p.d.f.) of the product $B^{-1}X$ of real 
and complex random matrices.

In order to motivate our discussion, before delving into random matrix 
theory, let us discuss a simpler problem. Thus, consider the random 
variables $x$ and $y$ drawn from the normal distribution 
\beq\label{normal}
G(x,y) = {1\over 2\pi \sigma^2} e^{-{x^2 + y^2\over 2\sigma^2}}\,.
\eeq
Define the random variable $z=\frac{x}{y}$. Obviously, its p.d.f. is 
independent of the width $\sigma$ of $G(x,y)$, and it is  a straightforward 
exercise to show that 
\beq\label{cauchy}
P(z) = {1\over \pi}{1\over 1+z^2}\,,
\eeq
i.e., the standard Cauchy distribution.

A slightly more interesting generalization of (\ref{normal}) is to consider 
the family of joint probability density (j.p.d.) functions of the 
form 
\beq\label{joint}
G(x,y) = f(x^2+y^2)\,,
\eeq
where $f(u)$ is a given appropriate p.d.f., subjected to the normalization 
condition 
\beq\label{normalizationf}
\int\limits_0^\infty\,f(u) du = {1\over \pi}\,.
\eeq
A straightforward calculation of the p.d.f. of $z=\frac{x}{y}$ leads again 
to (\ref{cauchy}). Thus, the random variable $z=\frac{x}{y}$ is distributed 
according to (\ref{cauchy}), {\em independently} of the function $f(u)$. 
In other words, (\ref{cauchy}) is a {\em universal} probability density 
function.\footnote{We can generalize (\ref{joint}) somewhat further, by 
considering circularly asymmetric distributions $G(x,y)=f(ax^2+by^2)$ (with 
$a,b >0$ of course, and the r.h.s. of (\ref{normalizationf}) changed to 
$\sqrt{ab}/\pi$), rendering (\ref{cauchy}) a Cauchy distribution of 
width $\sqrt{b/a}$, independently of the function $f(u)$.} $P(z)$ is 
universal, essentially, due to rotational 
invariance of (\ref{joint}). More generally, $P(z)$ must be independent, of 
course, of any common scale of the distribution functions of $x$ and $y$.\\

We will now show that an analog of this universal behavior exists in random
matrix theory. Our interest in this problem stems from the recent
application of random matrix theory made in \cite{bffs} to calculate 
the complexity of an analog computation process \cite{lp}, which solves linear 
programming problems.


\section{The universal probability distribution of the product $B^{-1}X$ of 
real random matrices}
\setcounter{equation}{0}

Consider a real $m\times (m+n)$ random matrix $A$ with entries 
$A_{i\alpha}\quad (i=1,\ldots m;~\alpha=1, \ldots m+n)$. We take the j.p.d. 
for the $m(m+n)$ entries of $A$ as
\beq\label{jpd}
G(A) = f(\rmtr A A^T ) = f\left(\sum_{i,\alpha} A_{i\alpha}^2\right)\,,
\eeq
with $f(u)$ a given appropriate p.d.f.. From\footnote{We use the ordinary 
Cartesian measure $dA = d^{m(m+n)}A = \prod_{i\alpha}dA_{i\alpha}$. 
Similarly, $dB = d^{m^2}B$ and $dX = d^{mn}X$ for the matrices $B$ and $X$
in (\ref{partition}) and (\ref{defpz}).} 
\beq\label{norma}
\int G(A)~dA ~ = 1
\eeq 
we see that $f(u)$ is subjected to the normalization condition 
\beq\label{normalizationff}
\int\limits_0^\infty\, u^{{m(m+n)\over 2}-1}f(u) du = {2\over S_{m(m+n)}}\,,
\eeq
where 
\beq\label{sphere}
S_d = {2\pi^{\frac{d}{2}}\over \Gamma \left(\frac{d}{2}\right)}
\eeq
is the surface area of the unit sphere embedded in $d$ dimensions.
This implies, in particular, that $f(u)$ must decay {\em faster} than 
$u^{-m(m+n)/2}$ as $u\rightarrow\infty$, and also, that if $f(u)$ blows up as
$u\rightarrow 0+$, its singularity must be {\em weaker} than 
$u^{-m(m+n)/2}$. In other words, $f(u)$ must be subjected to the asymptotic
behavior 
\beq\label{asymptoticsf}
u^{m(m+n)/2}\,f(u)\rightarrow 0
\eeq
both as $u\rightarrow 0$ and $u\rightarrow\infty$.\\
We now choose $m$ columns out of the $m+n$ columns of $A$, and pack them 
into an $m\times m$ matrix $B$ (with entries $B_{ij}$). Similarly, we pack the 
remaining $n$ columns of $A$ into an $m\times n$ 
matrix $X$ (with entries $X_{ip}$). This defines a partition 
\beq\label{partition} 
A\rightarrow (B,X)
\eeq
of the columns of $A$.

The index conventions throughout this paper are such that indices 
$$ i,j,\ldots\quad\quad {\rm range~over}\quad\quad  1,2,\ldots,m\,,$$
\beq\label{indices}
p,q,\ldots\quad\quad {\rm range~over}\quad\quad 1,2,\ldots, n\,,
\eeq
and $\alpha $ ranges over $1,2,\ldots, m+n\,.$

In this notation we have 
$\rmtr A A^T  = \sum_{i,j} B_{ij}^2 + \sum_{i,p} X_{ip}^2 = 
\rmtr B B^T  + \rmtr X X^T $, and thus (\ref{jpd}) reads 
\beq\label{bx}
G(B,X) = f(\rmtr B B^T  + \rmtr X X^T  )\,.
\eeq
We now define the random matrix $Z=B^{-1}X$. 
Our goal is to calculate the j.p.d. $P(Z)$ for the $mn$ entries of $Z$. 
$P(Z)$ is clearly independent of the particular partitioning 
(\ref{partition}) of $A$, since $G(B,X)$ is manifestly independent of 
that partitioning. The main result in this section is stated as follows: 
\begin{theorem}\label{realcase}
The j.p.d. for the $mn$ entries of the real random matrix $Z = B^{-1} X$ 
is independent of the function $f(u)$ and is given by the universal function 
\beq\label{result}
P(Z) = \frac{C}{[\det (1\!\! 1 + ZZ^T)]^{m+n\over 2}}\,,
\eeq
where $C$ is a normalization constant.
\end{theorem}
\begin{remark}\label{multivariate}
The probability density function (\ref{result}) is a special (spherically
symmetric) case of the so-called\footnote{Our notations in Remark 
(\ref{multivariate}) are slightly different from the notations used in
\cite{mvd}. In particular, we interchanged their $\Sigma$ and $\Omega$, and 
also denoted their $(T-M)^T$ by $Z-M$ here. Finally, we applied the identity
$\det (1\!\! 1 + AB) = \det (1\!\! 1 + BA)$ to arrive after all these 
interchanges from their equation (4.2.1) to (\ref{mvtdist}).} matrix variate 
$t$-distributions \cite{dickey,mvd}: The $m\times n$ random matrix $Z$ is 
said to have a matrix variate $t$-distribution with parameters 
$M,\Sigma,\Omega$ and $q$ 
(a fact we denote by $Z\sim T_{n,m}(q,M,\Sigma,\Omega)$) if its p.d.f. is 
given by 
\beq\label{mvtdist} 
D (\det\Sigma)^{-\frac{n}{2}}(\det\Omega)^{-\frac{m}{2}}
\left[\det\left(1\!\! 1_m + \Sigma^{-1}(Z-M)\Omega^{-1}(Z-M)^T\right)
\right]^{-{1\over 2}(m+n+q-1)}\,,
\eeq
where $M,\Sigma$ and $\Omega$ are fixed real matrices of dimensions 
$m\times n$, $m\times m$ and $n\times n$, respectively. $\Sigma$ and $\Omega$ 
are positive definite, and $q>0$. The normalization coefficient is 
\beq\label{tnormalization}
D={1\over\pi^{mn\over 2}}{\prod_{j=1}^n\Gamma\left({m+n+q - j\over 2}\right)
\over\prod_{j=1}^n \Gamma\left({n+q - j\over 2}\right)}\,.
\eeq
It arises in the theory of matrix variate distributions as the p.d.f. of a
random matrix which is the product of the inverse sqaure root of a certain
Wishart-distributed matrix and a matrix taken from a normal distribution,
and by shifting this product by $M$, as described in \cite{dickey,mvd}.
Our universal distribution (\ref{result}) corresponds to setting 
$M=0, \Sigma=1\!\! 1_m,\Omega=1\!\! 1_n$ and $q=1$ in (\ref{mvtdist}) and
(\ref{tnormalization}). 
\end{remark}
\begin{remark}
It would be interesting to distort the parent 
j.p.d. (\ref{jpd}) into a non-isotropic distribution and see if the 
generic matrix variate $t$-distribution (\ref{mvtdist}) arises as the 
corresponding universal probability distribution function in this case.
\end{remark}
To prove Theorem (\ref{realcase}), we need 
\begin{lemma}\label{clemma}
Given a function $f(u)$, subjected to (\ref{normalizationff}), the integral  
\beq\label{integral}
I=\int dB f(\rmtr\, B B^T)\, |\mbox{det}B|^n 
\eeq
converges, and is {\em independent} of the particular function $f(u)$. 
\end{lemma}
\begin{remark}
A qualitative and simple argument, showing the convergence of (\ref{integral}),
is that the measure $d\mu(B) = dB \,|\mbox{det}B|^n$ scales as $d\mu(tB) = 
t^{m(m+n)}d\mu(B)$, and thus has the same scaling property as 
$dA$ in (\ref{norma}), indicating that the integral 
(\ref{integral}) converges, in view of (\ref{asymptoticsf}). To see that 
I is independent of $f(u)$ one has to work harder.
\end{remark}
\begin{proof}
We would like first to integrate over the rotational degrees of freedom 
in $dB$. Any real $m\times m$ matrix $B$ may be decomposed as 
\cite{hua, rmt}
\begin{equation}\label{decomposition}
B = {\cal O}_1\Omega{\cal O}_2\,
\end{equation}
where ${\cal O}_{1,2}\in {\cal O}(m)$, the group of $m\times m$ orthogonal 
matrices, and $\Omega = {\rm Diag}(\omega_1,\ldots ,\omega_m)$, 
where $\omega_1,\ldots,\omega_m$ are the singular values of $B$. 
Under this decomposition we may write the measure $dB$
as \cite{hua, rmt}
\begin{equation}\label{measure}
dB = d\mu ({\cal O}_1) d\mu ({\cal O}_2) 
\prod_{i<j} | \omega_i^2-\omega_j^2 | d^m\omega\,,
\end{equation}
where $d\mu ({\cal O}_{1,2})$ are Haar measures over the 
$\cO (m)$ group manifold. The measure $dB$ is manifestly invariant under 
actions of the orthogonal group ${\cal O} (m)$
\begin{equation}\label{invariance}
dB=d(B{\cal O})=d({\cal O} 'B)\,,\quad\quad {\cal O}, {\cal O} ' \in 
{\cal O} (m)\,,
\end{equation}
as should have been expected to begin with.
\begin{remark}
Note that the decomposition (\ref{decomposition}) is not unique, since 
$\cO_1\cd$ and $\cd\cO_2$, with $\cd$ being any of the $2^m$ diagonal matrices
${\rm Diag}\,(\pm 1, \cdots, \pm 1)$, is an equally good pair of orthogonal 
matrices to be used in (\ref{decomposition}). Thus, as $\cO_1$ and $\cO_2$ 
sweep independently over the group $\cO (m)$, the measure (\ref{measure}) 
over counts $B$ matrices. This problem can be easily rectified by appropriately
normalizing the volume ${\cal V}_m = \int d\mu ({\cal O}_1) d\mu ({\cal O}_2)$.
One can show\footnote{One simple way to establish (\ref{volume}), is to 
calculate $\int dB\, \exp\,-\frac{1}{2}\rmtr B^T B = (2\pi)^{\frac{m^2}{2}} = 
{\cal V}_m\int\limits_{-\infty}^\infty d^m\omega\,
\prod_{i<j} | \omega_i^2-\omega_j^2 |\,\exp\,-\frac{1}{2}\sum_i \om_i^2$\,.
The last integral is a known Selberg type integral \cite{mehta}.}
that the correct normalization of the volume is 
\beq\label{volume}
{\cal V}_m = {\pi^{\frac{m(m+1)}{2}}\over 2^m\,\prod_{j=1}^m 
\Gamma\left(1+\frac{j}{2}\right)\Gamma\left(\frac{j}{2}\right)}\,.
\eeq
\end{remark}
Let us now turn to (\ref{integral}). The integrals over the orthogonal 
group in (\ref{integral}) clearly factor out, and we obtain 
\beq\label{singularvalues}
I={\cal V}_m\int\limits_{-\infty}^\infty \prod_{i=1}^m d\om_i~ \prod_{j<k} 
|\om_j^2-\om_k^2|\, 
\left(\prod_{i=1}^m \om_i\right)^n\,f\left(\sum_{i=1}^m \om_i^2\right)\,.
\eeq

Finally, we change the integration variables in (\ref{singularvalues}) to 
the polar coordinates associated with the $\om_i$. The angular part of that 
integral is fixed only by dimensionality and by the factor $\prod_{j<k} 
|\om_j^2-\om_k^2|\, \left(\prod_{i=1}^m \om_i\right)^n$, and is thus 
{\em independent} of the function $f(u)$.

To prove that $I<\infty$ we need only consider integration over the radius 
$r^2 = \sum_{i=1}^m \om_i^2$, since integration over the angles obviously 
produces a finite result. Using (\ref{normalizationff}), we find that 
the radial integral in question is 
$$\int\limits_0^\infty\,dr\,r^{m^2 + nm -1}\, f(r^2) = 
\frac{1}{2}\int\limits_0^\infty\, u^{{m(m+n)\over 2}-1}f(u) du  = 
{2\over S_{m(m+n)}}\,,$$ 
independently of $f(u)$. 
\end{proof}

We are ready now to prove Theorem (\ref{realcase}):

\begin{proof}
By definition\footnote{Our notation is such that $\delta(X)=\prod_{i=1}^m
\prod_{p=1}^n\delta (X_{ip}) = \prod_{p=1}^n\delta^{(m)}(X_p)$, 
$X_p$ being the $p$-th column of $X$.},  
\begin{eqnarray}
P(Z) & = & \int dB \,dX \, f( \rmtr \, B B^{T} + \rmtr \, X X^T \,) 
\delta (Z - B^{-1}X) 
\nonumber\\
&  = & \int dB \,dX \, f( \rmtr \, B B^{T}  + \rmtr \, X X^T  \,)
|\mbox{det}B|^n \, \delta (X - BZ) \,.
\label{defpz}
\end{eqnarray}
Integration over $X$ gives:
\begin{eqnarray}
\label{interim}
P(Z) & = & 
 \int dB \, f( \rmtr \, B B^{T}  + \rmtr \, B Z Z^T B^T) \,
|\mbox{det}B|^n \nonumber\\
 & = & \int dB \, f[ \rmtr \, B(1\!\! 1 + Z Z^T ) B^T] \, 
|\mbox{det}B|^n 
\end{eqnarray}
The $m\times m$ symmetric matrix $1\!\! 1 + Z Z^T $ 
can be diagonalized as $1\!\! 1 + Z Z^T = {\cal O} \Lambda {\cal O}^T$, 
where ${\cal O}$ is an orthogonal matrix, and $\Lambda = \mbox{Diag}
(\lambda_1,\ldots,\lambda_m)$ is the corresponding diagonal form. Obviously, 
all $\lambda_i\geq 1$, since $Z Z^T $ is positive definite.
Substituting this diagonal form into (\ref{interim}) we obtain 
\begin{eqnarray*}
P(Z) =  \int dB \, f( \rmtr \, B\cO \Lambda \cO^T B^T) \, 
|\mbox{det}B|^n \,.
\end{eqnarray*}

From the invariance of the determinant $|\mbox{det}B\cO| = |\mbox{det}B|$ and 
of the volume element $d(B\cO) = dB$ under orthogonal transformations we have:
\beq\label{interim1}
P(Z) = \int dB \, f( \rmtr \, B\Lambda B^T) \, |\mbox{det}B|^n \,.
\eeq
Let us now rescale $B$ as $\tilde B = B \sqrt{\Lambda}$. 
Thus, 
\beq\label{rescale}
\det\tilde B = \sqrt{\det\Lambda}\,\det B\,,\quad {\rm and}\quad 
d\tilde B = (\det\Lambda)^{m\over2}\,dB\,.
\eeq
Finally, substituting (\ref{rescale}) in (\ref{interim1}) we obtain 
\beqra\label{final}
P(Z) &=&  \int \frac{d\tilde{B}}{(\det\Lambda)^{m\over2}} \, 
f( \rmtr \, \tilde B \tilde{B}^T )
\,\left(\frac{|\mbox{det}\tilde{B}|}{(\det\Lambda)^{1\over2}}\right)^n
\nonumber\\
 & = & \frac{C}{(\det\Lambda)^{(m+n)/2}} = 
\frac{C}{[\det(1\!\! 1 + Z Z^T )]^{(m+n)/2}} \, ,
\eeqra
where $C$ is the normalization constant
\beq\label{c}
C=\int dB f(\rmtr\, B B^T)\, |\mbox{det}B|^n 
\eeq
rendering 
\beq\label{pz}
\int P(Z) dZ = 1\,.
\eeq
$C$ is nothing but the integral (\ref{integral}). Thus, according to 
Lemma (\ref{clemma}), $C<\infty$ and is also
independent of the function $f(u)$. 
\end{proof}
\begin{remark}\label{singz}
The j.p.d. $P(Z)$ in (\ref{result}) is manifestly a symmetric function 
only of the eigenvalues of $Z Z^T$, and thus, a symmetric function only 
of the singular values of $Z$. 
\end{remark}
\begin{remark}\label{independence}
From the normalization condition (\ref{pz}) we obtain an alternative 
expression for the normalization constant (\ref{c}) as
\beq\label{alternativec}
\frac{1}{C} = {1\over \int dB f(\rmtr\, B B^T)\, |\mbox{det}B|^n } = 
\int {dZ\over [\det(1\!\! 1 + Z Z^T )]^{(m+n)/2}}\,,
\eeq
which is manifestly independent of the particular function $f(u)$, 
in accordance with Lemma (\ref{clemma}). The integral over the matrix $Z$ 
can be reduced to a multiple integral of the Selberg type \cite{hua, mehta}
over the singular values of the matrix $Z$, which can be carried out explicitly:
\beq\label{selberg}
\int {dZ\over [\det(1\!\! 1 + Z Z^T )]^{(m+n)/2}} = \pi^{\frac{mn}{2}}\,
\prod_{j=1}^n\,{\Gamma\left(\frac{j}{2}\right)\over 
\Gamma\left(\frac{m+j}{2}\right)}\,.
\eeq
\end{remark}

For particular choices of the function $f(u)$, we can use (\ref{alternativec})
to derive explicit integration formulas. For example, the function 
\beq\label{gaussianf}
f(u) = {e^{-u}\over \pi^{m(m+n)/2}}
\eeq
(i.e., the entries $A_{i\alpha}$ in (\ref{jpd}) are i.i.d. according to a 
normal distribution of variance $1/2$) satisfies 
(\ref{normalizationff}). Thus, we obtain from (\ref{alternativec}) that  
\beq\label{nice}
\int  dB e^{-\rmtr\, B B^T}\, |\mbox{det}B|^n = 
\pi^{\frac{m^2}{2}}\,\prod_{j=1}^n\,{\Gamma\left(\frac{m+j}{2}\right)\over 
\Gamma\left(\frac{j}{2}\right)}\,.
\eeq
Note that the integral on the left-hand side of (\ref{nice}) can 
also be reduced to a multiple integral of the Selberg type (this time, over 
the singular values of $B$), which can be carried out explicitly. The result is 
$$\pi^{\frac{m^2}{2}}\,\prod_{j=1}^m\,{\Gamma\left(\frac{n+j}{2}\right)\over 
\Gamma\left(\frac{j}{2}\right)}\,.$$
Since this must coincide with (\ref{nice}), we obtain the identity 
\beq\label{amusing}
\prod_{j=1}^m\,{\Gamma\left(\frac{n+j}{2}\right)\over 
\Gamma\left(\frac{j}{2}\right)} = \prod_{j=1}^n\,{\Gamma\left(
\frac{m+j}{2}\right)\over \Gamma\left(\frac{j}{2}\right)}\,.
\eeq
\begin{example}\label{mcauchy}
For $n=1$, i.e., the case where $X$ and $Z$ are $m$ dimensional vectors, 
(\ref{result}) simplifies into the $m$ dimensional Cauchy distribution 
\beq\label{neq1}
P(Z) =  \frac{C}{(1 + Z^T Z )^{(m+1)/2}} \,.
\eeq
This is so because for $n=1$, the matrix $Z Z^T$ has $m-1$ eigenvalues 
equal to 0, that correspond to the $m-1$ dimensional subspace of vectors 
orthogonal to $Z$, and one eigenvalue equal to $Z^T Z$. Thus, 
$\det(1\!\! 1 + Z Z^T ) = 1 + Z^T Z$. Eq. (\ref{neq1}) then follows by 
substituting this determinant into (\ref{final}). 
\end{example}

\section{The universal probability distribution of the product $B^{-1}X$ of 
complex random matrices}
\setcounter{equation}{0}

The results of the previous section are readily generalized to complex 
random matrices. One has only to count the number of independent real 
integration variables correctly. In what follows we 
will use the notations defined in the previous section (unless specified
otherwise explicitly). Thus, consider a 
complex $m\times (m+n)$ random matrix $A$ with entries 
$A_{i\alpha}\quad (i=1,\ldots m;~\alpha=1, \ldots m+n)$. We take the j.p.d. 
for the $m(m+n)$ entries of $A$ as
\beq\label{jpdc}
G(A) = f(\rmtr A A^\dagger ) = f\left(\sum_{i,\alpha} |A_{i\alpha}|^2\right)\,,
\eeq
with $f(u)$ a given appropriate p.d.f.. From\footnote{We use the 
Cartesian measure $dA = d^{2m(m+n)}A = \prod_{i\alpha}d{\rm Re}\,A_{i\alpha}
~d{\rm Im} A_{i\alpha}$, with analogous definitions for $dB$ and $dX$ below.} 
\beq\label{normac}
\int G(A)~dA ~ = 1
\eeq 
we see that $f(u)$ is subjected to the normalization condition 
\beq\label{normalizationffc}
\int\limits_0^\infty\, u^{m(m+n)-1}f(u) du = {2\over S_{2m(m+n)}}\,.
\eeq
This implies that $f(u)$ must be subjected to the asymptotic
behavior 
\beq\label{asymptoticsfc}
u^{m(m+n)}\,f(u)\rightarrow 0
\eeq
both as $u\rightarrow 0$ and $u\rightarrow\infty$.\\
As in the previous section, we choose a partition 
\beq\label{partitionc} 
A\rightarrow (B,X)
\eeq
of the columns of $A$. Thus, $\rmtr A A^\dagger  = \sum_{i,j} |B_{ij}|^2 + 
\sum_{i,p} |X_{ip}|^2 = 
\rmtr B B^\dagger  + \rmtr X X^\dagger $, and thus (\ref{jpdc}) reads 
\beq\label{bxc}
G(B,X) = f(\rmtr B B^\dagger  + \rmtr X X^\dagger  )\,.
\eeq
We now define the random matrix $Z=B^{-1}X$. 
Our goal is to calculate the j.p.d. $P(Z)$ for the $mn$ entries of $Z$. 
The main result in this section is stated as follows: 
\begin{theorem}\label{complexcase}
The j.p.d. for the $mn$ entries of the complex random matrix $Z = B^{-1} X$ 
is independent of $f(u)$ and is given by the universal function 
\beq\label{resultc}
P(Z) = \frac{C}{[\det (1 + ZZ^\dagger)]^{m+n}}\,,
\eeq
where $C$ is the normalization constant (\ref{cc}).
\end{theorem}

\begin{proof}
The proof proceeds in a similar manner to the proof of 
Theorem (\ref{realcase}). The only important difference is that now 
$\delta(X)=\prod_{i=1}^m\prod_{p=1}^n\,\delta ({\rm Re}\,X_{ip})
\,\delta ({\rm Im}\,X_{ip}) = \prod_{p=1}^n\delta^{(2m)}(X_p)$, 
$X_p$ being the $p$-th column of $X$. One obtains 
\begin{eqnarray}\label{interimc}
P(Z) & = & \int dB \,dX \, f( \rmtr \, B B^{\dagger} + \rmtr \, X X^\dagger 
\,) \delta (Z - B^{-1}X) 
\nonumber\\
&  = & 
\int dB \, f[ \rmtr \, B(1\!\! 1 + Z Z^\dagger ) B^\dagger] \, 
|\mbox{det}B|^{2n} 
\end{eqnarray}
where we have integrated over $X$.

The $m\times m$ complex hermitean matrix $1\!\! 1 + Z Z^\dagger$ 
can be diagonalized as $1\!\! 1 + Z Z^\dagger = {\cal U} \Lambda 
{\cal U}^\dagger$, where ${\cal U}$ is a unitary matrix, and 
$\Lambda = \mbox{Diag}
(\lambda_1,\ldots,\lambda_m)$ is the corresponding diagonal form. Obviously, 
all $\lambda_i\geq 1$, since $Z Z^\dagger $ is positive definite.
Substituting this diagonal form into (\ref{interimc}) we obtain 
\beq\label{interim1c}
P(Z) =  \int dB \, f( \rmtr \, B\cu \Lambda \cu^\dagger B^\dagger) \, 
|\mbox{det}B|^{2n} = 
\int dB \, f( \rmtr \, B\Lambda B^\dagger) \, |\mbox{det}B|^{2n} \,,
\eeq
where we used the invariance of the determinant $|\mbox{det}B\cu| = 
|\mbox{det}B|$ and the invariance of the volume element $d(B\cu) = dB$ under 
unitary transformations. 

As in the previous section we now rescale $B$ as 
$\tilde B = B \sqrt{\Lambda}$. Thus, 
\beq\label{rescalec}
\det\tilde B = \sqrt{\det\Lambda}\,\det B\,,\quad {\rm and}\quad 
d\tilde B = (\det\Lambda)^{m}\,dB\,.
\eeq
Finally, substituting (\ref{rescalec}) in (\ref{interim1c}) we obtain 
\beqra\label{finalc}
P(Z) &=&  \int \frac{d\tilde{B}}{(\det\Lambda)^{m}} \, 
f( \rmtr \, \tilde B \tilde{B}^\dagger)
\,\left(\frac{|\mbox{det}\tilde{B}|}{(\det\Lambda)^{1\over2}}\right)^{2n}
\nonumber\\
 & = & \frac{C}{(\det\Lambda)^{m+n}} = 
\frac{C}{[\det(1\!\! 1 + Z Z^\dagger )]^{m+n}} \, ,
\eeqra
where $C$ is the normalization constant
\beq\label{cc}
C=\int dB f(\rmtr\, B B^\dagger)\, |\mbox{det}B|^{2n} 
\eeq
rendering 
\beq\label{pzc}
\int P(Z) dZ = 1\,.
\eeq
Finally, one can show, in a manner analogous to Lemma (\ref{clemma}), 
that $C<\infty$ and that it is independent of the particular function $f(u)$. 
\end{proof}

There are obvious analogs to the remarks made in the previous 
section, which follow from Theorem (\ref{complexcase}), which we will not 
write down explicitly. 

\section{More on the distribution of solutions of systems of linear equations 
with random coefficients: extension of a result due to Girko}
\setcounter{equation}{0}

The methods of the previous sections may be applied in studying the 
distribution of solutions of systems of linear equations 
with random coefficients. For concreteness, let us concentrate 
on real linear systems in real variables.

Consider a system of $m$ real linear equations 
\beq\label{system}
\sum_{\alpha=1}^{m+n}A_{i\alpha}\xi_\alpha = b_i\,,\quad i=1,\ldots m
\eeq
in the $m+n$ real variables $\xi_\alpha$. With no loss of generality, we will 
treat the first $m$ components of the vector $\xi$ as 
unknowns, and the remaining $n$ components of $\xi$ as given parameters. 
Thus, we split 
\beq\label{xisplit}
\xi = \left(\begin{array}{c} z\\u\end{array}\right)
\eeq
where $z$ is the vector of unknowns $z_i = \xi_i$ ($i=1,\ldots m$) and $u$ is 
the  vector of parameters $u_p = \xi_{p+m}$ ($p=1,\ldots n$). 
Similarly, we split the matrix of coefficients
$$ A = (B,X)\,, $$
where the $m\times m$ matrix $B$ (with entries $B_{ij}$) and the  
$m\times n$ matrix $X$ (with entries $X_{ip}$) were defined in (\ref{partition}).
Thus, we may rewrite (\ref{system}) explicitly as a system for the $z_i$:  
\beq\label{splitsystem}
B z = b - X u\,.
\eeq
If we consider an ensemble of systems (\ref{system}), in which $A$ and $b$ are 
drawn according to some probability law, the unknowns $z_i$ become random 
variables, which depend on the parameters $u_p$. Girko proved\cite{Girko} the 
following theorem for a particular family of such ensembles: 
\begin{theorem}\label{girko} (Girko)\\
If the random variables $A_{i\alpha}$ and $b_i$ 
($i=1,\ldots m;~\alpha=1, \ldots m+n $) are independent, identically 
distributed variables, having a stable distribution law with the
characteristic function 
\beq\label{characteristic}
g(t;\alpha, c) = e^{-c|t|^\alpha}\,,\quad 0<\alpha\leq 2;\, c>0\,,
\eeq
then the random variables $z_i$ ($i=1,\ldots m $) are identically distributed 
with the probability density function 
\beq\label{girkodensity}
p(\zeta; \alpha, \beta)  = \frac{2}{\beta}\int\limits_0^\infty\, r \rho \left(
\frac{r\zeta}{\beta};\alpha\right)\rho(r;\alpha)\, dr\,,
\eeq
where $\rho(r;\alpha)$ is the probability density of the postulated stable 
distribution, and 
\beq\label{beta}
\beta = \left( 1+ \sum_{p=1}^n |u_p|^\alpha\right)^{1\over\alpha}\,.
\eeq
The ratios $z_i/z_j$ ($i\neq j; i,j = 1,\ldots m $) have the density 
$p(r;\alpha,1)$. 
\end{theorem}
In the special case $\alpha = 2$ in (\ref{characteristic}), the random 
variables $A_{i\alpha}$ and $b_i$ are normally distributed. For this case
Girko obtained
\begin{corollary}\label{normalgirko}
If, under the conditions of Theorem (\ref{girko}), $\alpha =2$, then 
\beq\label{cauchygirko}
p(\zeta; 2, \beta)  = \frac{\beta}{\pi (\zeta^2 + \beta^2)}\,,
\eeq
i.e., $\zeta$ follows a Cauchy distribution of width $\beta$. 
\end{corollary}
When $\alpha=2$, the j.p.d. of the $A_{i\alpha}$ and the $b_i$ is  
\beq\label{jpdgirkonormal}
G(A,b) = (2\pi\si^2)^{\frac{m(m+n+1)}{2}}\, 
e^{-\frac{1}{2\si^2}\,(\rmtr A^T A + b^T b)}\,,
\eeq
which is a special case of the j.p.d.'s we have discussed in the previous 
sections. Thus, in the spirit of the discussion in the previous sections, 
we will study systems of linear equations (\ref{system}) with random 
coefficients $A$ and inhomogeneous terms $b$ with j.p.d.'s of the form
\beq\label{jpdgirko}
G(A, b) = f(\rmtr A^T A + b^T b)\,,
\eeq
with $f(u)$ a given appropriate p.d.f. subjected to the normalization condition
\beq\label{normalizationgirko}
\int\limits_0^\infty\, u^{{m(m+n+1)\over 2}-1}f(u) du = {2\over S_{m(m+n+1)}}\,.
\eeq

Our goal is to calculate the j.p.d. $P(z;u)$ for the $m$ unknowns $z_i$. 
We summarize our main result in this section as 
\begin{theorem}\label{universalgirko}
If the random variables $A_{i\alpha}$ and $b_i$ 
($i=1,\ldots m;~\alpha=1, \ldots m+n $) are distributed with a j.p.d. given
by (\ref{jpdgirko}), with $f(u)$ being any appropriate probability density 
function subjected to (\ref{normalizationgirko}), then 
the random variables $z_i$ ($i=1,\ldots m $) are distributed 
with the universal j.p.d. function 
\beq\label{universalgirkodensity}
P(z; u)  = C\, {\beta \over (\beta^2 + z^T z)^{\frac{m+1}{2}}}\,,
\eeq
independently of the function $f(u)$, where 
\beq\label{beta1}
\beta = \sqrt{1 + u^T u}
\eeq
and $C$ is a normalization constant given by 
\beq\label{cgirko}
C = \int\,dA\, |{\rm det} B|\,f(\rmtr A^T A)\,.
\eeq
\end{theorem}
\begin{remark}
Note that Theorem (\ref{universalgirko}) generalizes the case $\alpha =2$ of 
Girko's result, Theorem (\ref{girko}), from the particular 
$f(u)\sim e^{-u}$ to a whole class of probability densities $f(u)$, and 
moreover, it determines for this class of distributions the (universal) 
j.p.d. of the $z_i$'s, and not only the distribution of a single component. 
Thus, it is an interesting question whether Girko's result could be 
generalized also to other ensembles of systems of linear equations as well. 
\end{remark}
\begin{proof}
By definition, from (\ref{splitsystem}),
\begin{eqnarray}
P(z;u) & = & \int dA \,db \, f( \rmtr \, A^T A + b^T b \,) 
\delta (z - B^{-1}(b- X u)) 
\nonumber\\
& = & \int dA \,db \, f( \rmtr \, A^T A + b^T b \,) 
|\mbox{det}B| \, \delta (Bz + Xu -b)
\nonumber\\
&  = & \int dA \,f[ \rmtr \, A^T A + (A\xi)^T (A\xi)\,]
|\mbox{det}B|\,,
\label{defpzgirko}
\end{eqnarray}
where in the last step we integrated over the $m$ dimensional vector $b$ and 
used $Bz + Xu = A\xi$. 

The last expression in (\ref{defpzgirko}) is manifestly 
invariant under $\cO (m)\times \cO (n) $ orthogonal transformations 
\beq\label{omn}
P(\cO_1\, z; \cO_2\, u) = P(z;u)\,,\quad\quad
\cO_1\in\cO(m)\,,\cO_2\in\cO(n),
\eeq
due to the invariance of the measure $dA = dB\,dX = d(B\cO_1)\,d(X\cO_2)$, 
the invariance of the determinant $|{\rm det}\,B| = |{\rm det}\, B\cO_1|$, 
and the invariance of the trace $\rmtr \, A^T A = \rmtr \, (B\cO_1)^T (B\cO_1)
+ \rmtr \, (X\cO_2)^T (X\cO_2).$ With this symmetry at our disposal, we may
simplify the calculation of $P(z;u)$ by rotating the vectors $z$ and $u$ into 
fixed convenient directions, e.g., into the directions in which only $z_1$ and 
$u_1$ do not vanish:
\beq\label{rotated}
z_i^{(0)} = z_0\,\delta_{i1}\,,\quad 
u_p^{(0)} = u_0\,\delta_{p1}\,.
\eeq
with 
\beq\label{z0u0}
z_0 = (z^T z)^{\frac{1}{2}}\,,\quad 
u_0 = (u^T u)^{\frac{1}{2}}\,.
\eeq
Thus, we obtain
\beqra\label{p0}
P(z;u) & = & \int dA \,f[ \rmtr \, A^T A + (Bz^{(0)} + Xu^{(0)})^T 
(Bz^{(0)} + Xu^{(0)})\,]\,|\mbox{det}B|
\nonumber\\
& = & \int dB\,dX \,f(S\,)\,
|\mbox{det}B|\,,
\eeqra
where
\beq\label{s}
S=\sum_{j=2}^m\,B_j^T B_j + 
\sum_{p=2}^n\,X_p^T X_p  
+ (1+z_0^2)\,B_1^T B_1  + 2u_0z_0\, B_1^T\,X_1 + (1+u_0^2)\,X_1^T X_1\,,
\eeq
in which $B_i$ is the $i$-th column of $B$, and $X_p$ is the $p$-th column of 
$X$. 

The bilinear form involving $B_1$ and $X_1$ in (\ref{s}) may be diagonalized as 
\beq\label{bilinear}
(1+\xi^T \xi) \left( {z_0 B_1 + u_0 X_1\over \sqrt{\xi^T \xi}}\right)^2
 + \left( {u_0 B_1 -z_0 X_1\over \sqrt{\xi^T \xi}}\right)^2\,,
\eeq
where we have used $z_0^2 + u_0^2 = \xi^T \xi$. We now perform a rotation
in the $B_1-X_1$ plane, followed by a scale transformation of the first term 
in (\ref{bilinear}), thus defining
\beq\label{b'x'}
B_1' = (1+\xi^T \xi)^{\frac{1}{2}}\, {z_0 B_1 + u_0 X_1\over \sqrt{\xi^T \xi}}\,,
\quad\quad
X_1' = {u_0 B_1 -z_0 X_1\over \sqrt{\xi^T \xi}}\,,
\eeq
such that $d^m B_1'\,d^m X_1' = (1+\xi^T \xi)^{\frac{m}{2}}\,d^m B_1\,d^m X_1$. 
We will also need the inverse transformation for $B_1$
\beq\label{inverseb}
B_1 (B_1', X_1') = {1\over \sqrt{\xi^T \xi}}\,\left({z_0B_1'\over 
(1+\xi^T \xi)^{\frac{1}{2}}}+ u_0X_1'\right)
\eeq
in order to express the matrix $B$ in terms of the primed column vectors:
\beq\label{btilde}
\tilde B = \left( B_1 (B_1', X_1'), B_2, \ldots, B_m\right)\,.
\eeq
Thus, using (\ref{bilinear}) - (\ref{btilde}) and the trivial fact that 
$dA = \prod_{i=1}^m\, d^mB_i\,\prod_{i=p}^n\, d^mX_p$, we obtain 
\beq\label{p00}
P(z;u) =  {1\over (1+\xi^T \xi)^{\frac{m}{2}}}\,\int dA \,
f(\rmtr \, A^T A\,)\,|\mbox{det}\tilde B|\,,
\eeq
where we have removed the primes from the integration variables. 
We are not done yet, since $\tilde B_1$, the first column of $\tilde B$, 
depends on $z_0$ and $u_0$. To rectify this 
problem, we note from 
(\ref{inverseb}) that 
\beq\label{b1tilde}
\tilde B_1 (B_1, X_1) = \sqrt{{1+u_0^2\over 1+\xi^T \xi}}\,\left( B_1\,
\mbox{cos}\,\theta + X_1\,\mbox{sin}\,\theta\,\right)\,,\quad\quad 
\mbox{cos}\,\theta = {z_0\over \sqrt{\xi^T \xi\,(1+ u_0^2)}}\,.
\eeq
Thus, performing one final rotation by an angle $\theta$ in the $B_1-X_1$ plane,
which leaves, of course, $dA$ and $\rmtr A^T A$ invariant, we see that in 
terms of the rotated columns $|\mbox{det}\tilde B| = 
\sqrt{{1+u_0^2\over 1+\xi^T \xi}}\,|\mbox{det}B|$, and thus, finally,
we obtain that
\beq\label{finalgirko}
P(z;u) =  {\sqrt{1 + u_0^2}\over (1+\xi^T \xi)^{\frac{m+1}{2}}}\,\int dA \,
f(\rmtr \, A^T A\,)\,|\mbox{det} B|\,,
\eeq
which coincides with (\ref{universalgirkodensity}), due to (\ref{beta1}) and 
(\ref{cgirko}). 
\end{proof}
\begin{remark}
The fact that the integral 
$C = \int\,dA\, |{\rm det} B|\,f(\rmtr A^T A)$ is convergent and independent
of the function $f(u)$ can be proved by decomposing $A$ into its singular
values \cite{hua, rmt}, essentially in a manner similar to our proof of Lemma 
(\ref{clemma}), but with slight modifications in (\ref{decomposition}) and 
(\ref{measure}) due to the fact that $A$ is a rectangular matrix rather
than a square matrix. We shall not get into these technicalities here, 
which the
reader may find in \cite{hua, rmt}. Note, however, that $C$ may be determined
from the normalization of $P(z;u)$: 
$$\int\,d^m z P(z;u) = 1= C\beta\,\int\, 
{d^m z\over (\beta^2 + z^T z)^{\frac{m+1}{2}}} = {C\pi^{\frac{m+1}{2}}\over 
\Gamma\left(\frac{m+1}{2}\right)}\,.$$
Thus, 
\beq\label{cvalue}
C = {\Gamma\left(\frac{m+1}{2}\right)\over \pi^{\frac{m+1}{2}}} = 
\frac{2}{S_{m+1}}\,.
\eeq
\end{remark}
\begin{remark}
We note that for $n=0$ (i.e., when $A=B $) and $u=0$, (\ref{splitsystem})
degenerates into $B z = b$, which is precisely the case $n=1$ (and $z\equiv Z
$) in the conditions for Theorem (\ref{realcase}), which we analyzed in 
Example (\ref{mcauchy}). Thus, (\ref{universalgirkodensity}), evaluated 
at $n=1$ and $u=0$ must coincide with (\ref{neq1}), as one can easily 
check it does. 
\end{remark}

Since Theorem (\ref{universalgirko}) states the explicit form 
(\ref{universalgirkodensity}) of $P(z;u)$, we can now use it to derive, e.g.,  
the probability density of the distribution of a single component $z_i$ and 
that of the ratio of two different components, mentioned in Girko's Theorem 
(\ref{girko}): 

\begin{corollary}\label{1ptfunction}
The $m$ components $z_i$ are identically distributed, with the probability 
density of any one of the components $z_i = \zeta$ given by 
(\ref{cauchygirko}) of Corollary (\ref{normalgirko}).
\end{corollary}
\begin{proof}
That the $z_i$ are identically distributed is an immediate consequence of 
the rotational invariance of $P(z;u)$ in (\ref{universalgirkodensity}). 
The proof is completed by performing the necessary integrals:
\beqra\label{cauchygirko1}
p(\zeta; \beta)  &=& \int\, dz_2\ldots dz_m\, P(z;u)_{|_{z_1=\zeta}} 
= C\beta\, \int\, {dz_2\ldots dz_m \over (\beta^2 + \zeta^2 + 
\sum_{i=2}^m z_i^2)^{\frac{m+1}{2}}}\nonumber\\{}\nonumber\\
&=& {C\pi^{\frac{m-1}{2}}\over 
\Gamma\left(\frac{m+1}{2}\right)}\,{\beta\over \beta^2 + \zeta^2}\,.
\eeqra
Thus, from (\ref{cvalue}) we obtain the desired result that 
$p(\zeta; \beta) = {\beta\over \pi(\beta^2 + \zeta^2)}\,.$
This result, should have been anticipated, since the universal
formula (\ref{universalgirkodensity}) holds, in particular, for the 
the Gaussian distribution (\ref{jpdgirkonormal}). 
\end{proof}

Finally, we have 
\begin{corollary}\label{1ratios}
The ratios $z_i/z_j$ ($i\neq j; i,j = 1,\ldots m $) have the density 
$P(r) =  p(r;1)$. 
\end{corollary}
\begin{proof}
The ratio $z_i/z_j$ is dimensionless, and thus its distribution cannot
depend on the width $\beta$, which is the only dimensionful quantity in 
(\ref{universalgirkodensity}). The proof amounts to performing the 
necessary integrals, e.g., for the random  variable $z_1/z_2$:
\beqra\label{11ratio}
P(r)  &=& \int\, d^m z\, P(z;u) \delta \left( r- \frac{z_1}{z_2}\right)
= C\beta\,\int\, {|z_2|dz_2\ldots dz_m\over 
[\beta^2 + (r^2+1)z_2^2 + \sum_{i=3}^m z_i^2]^{\frac{m+1}{2}}}
\nonumber\\{}\nonumber\\
&=& {2CS_{m-2}\over r^2 +1}\, 
{\Gamma\left(\frac{m-2}{2}\right)\,\Gamma\left(\frac{3}{2}\right)\over 
2\Gamma\left(\frac{m+1}{2}\right)} = {1\over \pi(r^2 + 1)} = p(r;1)\,.
\eeqra
\end{proof}

\pagebreak

\section*{Acknowledgments}
I thank Asa Ben-Hur, Shmuel Fishman and Hava Siegelmann
for a stimulating collaboration on complexity in analog computation, 
which inspired me to consider the problems studied in this paper. I also 
thank Ofer Zeitouni for making some comments. Shortly after the first version 
of this work appeared, its connection with the theory 
of matrix variate distributions was pointed out to me by Peter Forrester. 
I am indebted to him for enlighting me on that connection and for suggesting 
references \cite{dickey,mvd}. This research was supported in 
part by the Israeli Science Foundation.

\end{document}